\documentclass[12pt]{amsart}
\usepackage{amssymb}

\usepackage{euscript}

\let\cal\mathcal

\makeatletter \@mparswitchfalse \makeatother

\newtheorem{theorem}{Theorem}
\newtheorem{lemma}{Lemma}
\newtheorem{corollary}{Corollary}

\newtheorem{remark}{Remark}
\newtheorem{remarks}{Remarks}

\newtheorem{definition}{Definition}

\newcommand{\pprec}{\prec\!\prec}

\def\eqref#1{(\ref{eq#1})}
\def\eqlabel#1{\label{eq#1}}

\def\N{\mathbb N}

\def\R{\mathbb R}
\def\M{\cal{M}}
\let\<\langle

\def\T{\tau}

\numberwithin{equation}{section}

\begin{document}


\let\\\cr
\let\union\bigcup
\let\inter\bigcap
\def\supp{\operatorname{supp}}
\def\Im{\operatorname{Im}}
\def\dim{\operatorname{dim}}
\def\Span{\operatorname{span}}
\def\cord{\operatorname{cord}}
\def\Re{\operatorname{Re}}
\def\sqn{\operatorname{sqn}}
\def\log{\operatorname{log}}
\def\max{\operatorname{max}}
\def\min{\operatorname{min}}
\let\emptyset\varnothing
\def\exp{\operatorname{exp}}

\title[Hilbert Transform]
{Hilbert Transform Associated with Finite Maximal Subdiagonal Algebras}
\author{Narcisse Randrianantoanina}
\address{Department of Mathematics and Statistics, Miami University, Oxford,
 OH 45056}
\email{randrin@muohio.edu}

\begin{abstract}
Let $\M$ be a von Neumann algebra with a faithful normal trace $\T$,
and let $H^\infty$ be a finite, maximal, subdiagonal algebra of $\M$.
Fundamental theorems on conjugate functions for weak$^*$\!-Dirichlet
algebras are shown to be valid for non-commutative $H^\infty$. In
particular the Hilbert transform is shown to be a bounded linear map
from $L^p(\M, \T)$ into $L^p(\M, \T)$ for $1 < p < \infty$, and to be a
continuous map from  $L^1(\M,\T)$ into $L^{1, \infty}(\M,\T)$. We also
obtain that if a 
 positive operator $a$ is such that $a\log^+a \in L^1(\M,\T)$
then its conjugate belongs to $L^1(\M,\T)$.
\end{abstract}

\maketitle

\section{introduction}
The theory of conjugate functions has been a strong motivating force behind
various aspects of harmonic analysis and abstract analytic function spaces.
This theory which was originally developed for functions in the circle group
$\mathbb{T}$ has found many generalizations to  more abstract  settings
 such as
Dirichlet algebras in \cite{DEV} and weak*-Dirichlet algebras in \cite{HR}.
Results from this theory
 have been proven to be very fruitful for studying Banach space
properties of the Hardy spaces (and their relatives) associated with the
 algebra involved (see for instance \cite{BO6} and \cite{LAN}).

Let $\M$ be a von Neumann algebra with a faithful, normal finite trace $\T$.
 Arveson introduced in \cite{AV}, as non-commutative
analogues of weak$^*$\!-Dirichlet algebras, the notion of finite, maximal
subdiagonal algebras of $\M$ (see definition below).
 Subsequently several authors
studied the (non-commutative) $H^p$-spaces associated with  such algebras
(\cite{KATO}, \cite{MAR}, \cite{MMS}, \cite{SAI2}, \cite{SAI1}).
 In \cite{MAR},  the notion of
harmonic conjugates  was introduced
for maximal subdiagonal algebras
 generalizing the notion of conjugate functions for
weak*-Dirichlet algebras and it was proved that the operation
of conjugation is bounded in $L^2(\M, \T)$.

The main objective of this paper  is to
combine the spirit of \cite{HR} with that of
\cite{MAR} to get a more constructive
 definition of conjugate operators for the setting of non-commutative
maximal subdiagonal algebras;
and  to study different properties of conjugations
 for these non-commutative settings.  We prove that most fundamental
theorems on conjugate operation on Hardy spaces associated with
weak$^*$ \!-Dirichlet (see \cite{DEV} and \cite{HR})
 remain valid for Hardy spaces associated with finite subdiagonal
algebras. In particular, we show that the conjugation operator is a
bounded map from $L^p(\M,\T)$ into $L^p(\M, \T)$ for $1< p < \infty$,
 and from $L^1(\M, \T)$ into $L^{1, \infty}(\M, \T)$.
 We conclude that, as in commutative case,
 (non-commutative) $H^p$ is a complemented subspace of  $L^p(\M, \T)$
for $1<p<\infty$.

We refer to \cite{N}, \cite{SEG} and  \cite{TAK} for general information
 concerning
von Neumann algebras as well as basic notions of non-commutative
integration,  to \cite{D1} and \cite{LT} for Banach space theory
and to \cite{HEL} and \cite{ZYG} for basic definitions from
harmonic analysis. %
\section{Definitions and preliminary results}

Throughout, $H$ will denote a Hilbert space and $\M \subseteq
\cal{L}(H)$ a von Neumann algebra with a normal, faithful finite trace
$\T$. A closed densely defined operator $a$ in $H$ is said to be {\em
affiliated with} $\M$ if $u^* au = a$ for all unitary $u$ in the
commutant $\M'$ of $\M$. If $a$ is a densely defined self-adjoint
operator on $H$, and if $a = \int^\infty_{- \infty} s d e^a_s$ is its
spectral decomposition, then for any Borel subset $B \subseteq \R$,
we denote by $\chi_B(a)$ the corresponding spectral projection
$\int^\infty_{- \infty} \chi_B(s) d e^a_s$. A closed densely defined
operator on $H$ affiliated with $\M$ is said to be {\em
$\T$-measurable} if there exists a number $s \geq 0$ such
that $\T(\chi_{(s, \infty)} (|a|)) < \infty$.

The set of all $\T$-measurable operators will be denoted by
$\overline{\M}$. The set $\overline{\M}$ is a $*$\!-algebra with respect to the
strong sum, the strong product, and the adjoint operation \cite{N}. %
For $x \in \overline{\M}$, the generalized singular value function $\mu (x)$
of $x$ is defined by
$$
\mu_t(x) = \inf \{ s \geq 0: \T(\chi_{(s, \infty)} (|x|)) \leq t \},
\quad \text{ for } t \geq 0.
$$
The function $t \to \mu_t(x)$ from $(0, \T(I))$ to $[0, \infty)$
is right continuous, non-increasing and is the inverse of the
distribution function $\lambda (x)$, where $\lambda_s(x) = \T(\chi_{(s,
\infty)}(|x|))$, for $s \geq 0$. For a complete study of $\mu(.)$ and
$\lambda(.)$,  we refer to \cite{FK}.

\begin{definition}
Let $E$ be an order continuous rearrangement invariant (quasi-)
Banach function space on $(0, \T(I))$. We define the symmetric space
$E(\M, \T)$ of measurable operators by setting:
\begin{align*}
E(\M, \T) &= \{ x \in
\overline{\M}\quad ; \quad \mu(x) \in E \} \quad \text{and} \\
\|x\|_{E(\M,\T)} &= \| \mu(x)\|_E, \text{ for } x \in E(\M,\T).
\end{align*}
\end{definition}

It is well known that $E(\M, \T)$ is a Banach space
(resp. quasi-Banach space) if $E$ is a Banach
space (resp. quasi-Banach space),
 and that if $E = L^p(0, \T(I))$, for $0 < p < \infty$, then
$E(\M, \T)$ coincides with the usual non-commutative $L^p$\!-space
associated with $(\M, \T)$. We refer to \cite{CS}, \cite{DDP1} and \cite{X}
 for more detailed discussions about these
spaces. For simplicity we will always assume that the trace $\tau$ is
normalized. %

The following definition isolates the main topic of this paper.

\begin{definition}
Let $H^\infty$ be a weak$^*$\!-closed unital subalgebra of $\M$ and let
$\Phi$ be a faithful, normal expectation from $\M$ onto the diagonal $D
= H^\infty \cap (H^\infty)^*$, where $(H^\infty)^* = \{ x^*,\  x \in
H^\infty \} $. Then $H^\infty$ is called a finite, maximal,
subdiagonal algebra in $\M$ with respect to $\Phi $ and $\T$ if:
\begin{itemize}
\item[(1)] $H^\infty + (H^\infty)^*$ is weak$^*$\!-dense in $\M$;
\item[(2)] $\Phi (ab) = \Phi(a) \Phi(b)$ for all $a,b \in H^\infty$;
\item[(3)] $H^\infty$ is maximal among those subalgebras satisfying (1)
and (2);
\item[(4)] $\T \circ \Phi = \T$.
\end{itemize}
\end{definition}

For $0 < p < \infty$, the closure of $H^\infty$ in $L^p(\M, \T)$ is
denoted by $H^p(\M,\T)$ (or simply $H^p$) and is called the Hardy space
 associated with the subdiagonal algebra $H^\infty$. Similarly,
  the closure of $H^\infty_0 = \{ x \in H^\infty;
\Phi(x) = 0 \}$ is denoted by $H^p_0$.

Note that $\Phi$ extends to $L^2(\M, \T)$ and this extension is an
orthogonal projection from $L^2(\M, \T)$ onto $[D]_2$, the closure of
$D$ in $L^2(\M, \T)$. Similarly, since $ \| \Phi(x)\|_1 \leq \|x\|_1$
for every $x \in \M$, the operator $\Phi$ extends
uniquely to a projection of norm one from  $L^1(\M, \T)$ onto
$[D]_1$, the closure of $D$ in $L^1(\M, \T)$.

\section{Harmonic conjugates and Hilbert transform}

Let $\cal{A} = H^\infty + (H^\infty)^*$. Since $\cal{A}$ is
weak$^*$\!-dense in $\M$, it
is norm dense in $L^p(\M, \T)$, where $ 1 \leq p < \infty$.

Note that $H^\infty$ and $(H^\infty_0)^*$ are orthogonal in $L^2(\M,\T)$.
This fact implies that $L^2(\M,\T)=H^2 \oplus (H^2_0)^*$, and hence
that $L^2(\M,\T)= H^2_0 \oplus (H^2_0)^* \oplus [D]_2$.

Let $a \in \cal{A}$. Then $a$ can be written as  $a_1 + a_2^* + d$ where
$a_1$ and   $a_2$ belong to $ H^\infty_0$ and $d \in D$.
 In fact, $ a = b_1 + b_2^*$
with $b_1, b_2 \in H^\infty $ and set $d = \Phi(b_1) + \Phi(b_2^*) \in
D$ and $a_i = b_i - \Phi(b_i)$, for $i = 1,2$.

Since $H^\infty_0$ and $(H^\infty_0)^*$ are orthogonal subsets of
$L^2(\M, \T)$, this decomposition is unique.

 For $u = u_1 + u_2^* + d$ in
$\cal{A}$ , we define $\tilde{u} = i u_2^* - i u_1$. Then $\tilde{u} \in \M$
and $u + i \tilde{u} = 2 u_1 + d \in H^\infty$.
The operator $\tilde{u}$ will be called the {\bf conjugate} of $u$ and the
{\bf Hilbert transform (conjugation operator)} $\cal{H}$ is defined on
 $\cal{A}$ as
follows: $\cal{H}: \cal{A}\to \M \, (u \to \cal{H}(u)=\tilde{u})$.

 Our main
objective is to study the Hilbert transform  as linear operator between
 non-commutative $L^p$-spaces.
 In particular we will
 extend   $\cal{H}$ to $L^p(\M, \T)$ for $1 \leq p <
\infty$. It should be noted that if $\M$ is commutative, then the above
definition coincides with the the definition of conjugate functions for
weak$^*$\!-Dirichlet algebras studied in \cite{HR}.

\begin{remark}
\begin{itemize}

\item[(i)] If $u = u^*$, then the uniqueness of the decomposition
implies that $u_1 = u_2$ and $d=d^*$.
 Therefore  if $u=u^*$ then $\tilde{u} = \tilde{u}^*$.
\item[(ii)] For $u = u_1 + u_2^* + d \in \cal{A}$ and
$\tilde{u} = i(u_2^* - u_1)$, the above observation implies that
$u_2^* \perp u_1 $ in $L^2(\M,\T)$, so
$$
\| \tilde{u} \|^2_2 = \| u_2^* - u_1 \|^2_2 = \| u_2^* \|^2_2 + \| u_1
\|^2_2,
$$
and since $L^2(\M,\T) = H^2_0 \oplus (H^2_0)^* \oplus [D]_2$ we get,
 $$\| u\|^2_2 = \| u_1 \|^2_2 + \| u_2^* \|^2_2 + \| d \|^2_2,$$
 which implies
that $\| \tilde{u} \|_2 \leq \| u \|_2$.
\end{itemize}
\end{remark}
As a consequence of (ii), we get the following theorem:
\begin{theorem}
There is a unique continuous linear map $\cal{H}$ from
$L^2(\M, \T)$ into $L^2(\M,\T)$ that coincides with  $\cal{H}$ in $\cal{A}$.
This map is of norm 1, and if $u \in L^2(\M, \T)$ then $ u + i\cal{H}(u)
\in H^2$.
\end{theorem}
We remark that Marsalli has recently proved a version of Theorem~1 (see
\cite{MAR} Corollary~10):
 he showed that the conjugation operator is bounded in $L^2(\M,\T)$
with bound less than or equal to $\sqrt{2}$.

\smallskip
Our next result is an extension of Theorem~1 from $p=2$ to all $p$ with
$1<p<\infty$.

\begin{theorem}
For each $1 < p < \infty$, there is a unique continuous linear
extension of $\cal{H}$ (which is also denoted by $\cal{H}$)
from $L^p(\M,\T)$ into $ L^p(\M,\T)$ with the
property that $f + i \tilde{f} \in H^p$ for all $f \in L^p(\M,\T)$ and
$\cal{H}(f)=\tilde{f}$.
Moreover there is a constant $C_p$ such that
$$
\| \tilde{f}\|_p \leq C_p  \| f \|_p
\quad \text{ for all } f \in L^p(\M,\T).
$$
\end{theorem}

The following elementary lemma will be used in the sequel; we
will include its  proof for completeness.

\begin{lemma}
Let $m \in \N$ and $a_1,a_2, \ldots, a_m \in \overline{\M}$. If
$\frac{1}{p_1} + \frac{1}{p_2} + \ldots + \frac{1}{p_m}=1$, and $a_j \in
L^{p_j}(\M,\T)$ for each $ j \leq m$, then
$$
|\T(a_1 a_2 \ldots a_m)| \leq \Pi^m_{j=1} \| a_j \|_{p_j}.
$$
\end{lemma}

\begin{proof}
 Recall that, for $a, b \in \overline{\M}$, the operator  $a$ is said to be
submajorided by $b$ and write $a \pprec b$ if
$$
\int^\alpha_0 \mu_t (a) dt \leq \int^\alpha_0 \mu_t (b) dt, \quad
\text{ for all } \alpha \geq 0.
$$
The lemma will  be proved inductively on $m \in \N$:

 For $m=2$, it is the usual
H\"older's  inequality.

Let $\frac{1}{p} + \frac{1}{q} = \frac{1}{r}$ and $b, c \in
\overline{\M}$. Then $\| bc \|_r \leq \| b \|_p \cdot \| c \|_q$: this is
a consequence of the fact \cite[Theorem 4.2\,(iii)]{FK} that
$
\mu_{(.)}^l (bc) \pprec \mu_{(.)}^l(b) \mu_{(.)}^l (c) \quad \text{
for all }l \in \N$.
So
$$
\| bc \|^r_r = \int \mu^r_t (bc)dt \leq \int \mu^r_t (b) \mu^r_t (c)
dt;
$$
then apply the usual H\"older's inequality for functions.

 Now assume
that the lemma is valid for $m=1, 2, \ldots, k$. Let $a_1, a_2,
\ldots, a_k, a_{k+1} \in \overline{\M}$ and $\frac{1}{p_1} +
\frac{1}{p_2} + \cdots + \frac{1}{p_k} + \frac{1}{p_{k+1}}
= 1$. Choose $q$ such that
$\frac{1}{q}=\frac{1}{p_k}+\frac{1}{p_{k+1}}$.  Then
\begin{align*}
\T(a_1\ldots a_{k-1} \cdot(a_k a_{k+1})) & \leq
\Pi^{k-1}_{i=1} \| a_j \|_{p_j} \cdot \| a_k
a_{k+1}\|_q \\
& \leq \Pi^{k+1}_{i=1} \| a_j \|_{p_j}.
\end{align*}
The proof is complete.
\end{proof}

\noindent
{\bf Proof of Theorem~2.}

Our proof follows Devinatz's argument (\cite{DEV}) for Dirichlet algebras,
but at number of points, certain
non-trivial adjustments have to be made to fit the non-commutative setting.

Let $u \in \cal{A}$ be  nonzero and self-adjoint; $\tilde{u}$ is self-adjoint.
Let $g=u+i \tilde{u} \in H^\infty$. Since $u = u^*$, it is of the form
$u = a+a^* + d$, where $a \in  H^\infty_0$ and $d=d^* \in D$. Recall that
$\tilde{ u}= i(a^* - a)$ so $g = 2a + d \in H^\infty$. We get that
$$
\Phi(g^{2k}) = \Phi ((2a+d)^{2 k}) = [2 \Phi(a) +
\Phi(d)]^{2 k } = \Phi (d)^{2 k}.
$$
So $\Phi((u+i \tilde{u})^{2 k}) = \Phi(d)^{2 k}$ and %
taking the adjoint, $\Phi((u - i \tilde{u})^{2 k}) = \Phi (d)^{2
k}$. Adding these two equalities, we get
\begin{equation}
\eqlabel{3}
\Phi [ ( u + i \tilde{u})^{2 k} + (u - i \tilde{u})^{2 k}] =
2 \Phi(d)^{2 k}.
\end{equation}
Now we will expand the operators $(u + i \tilde{u})^{2 k}$ and $(u - i
\tilde{u})^{2 k}$. Note that $u$ and $\tilde{u}$  do not
necessarily commute.

For $2\leq m\leq 2k$, let
 $\cal{S}_m = \{(r_1, r_2, \ldots , r_{m}) \in \{ 1, \ldots, 2k-1\}^{m };
\quad \sum^{m }_{j=1} r_j = 2k\}$
 and set $\cal{S}=\cup_{2\leq m \leq 2k} \cal{S}_m$.
For a finite sequence of integers $r=(r_1, r_2, \ldots, r_m)$,
we set $s(r)=\sum_{j=1}^{[m/2]} r_{2j}$. Then
\begin{align*}
(u + i \tilde{u})^{2k} &= u^{2k} +(i\tilde{u})^{2k} +
 \sum_{(r_1,  \ldots, r_{m})\in
\cal{S}} (u^{r_1} (i\tilde{u})^{r_2}\ldots) +
  ((i\tilde{u})^{r_1}u^{r_2} \ldots) \\
&= u^{2k} +(i)^{2k}\tilde{u}^{2k} +
\sum_{(r_1, \ldots, r_{m}) \in \cal{S}}(i)^{s(r)}
(u^{r_1} \tilde{u}^{r_2}\ldots)  +
 (i)^{2k - s(r)} (\tilde{u}^{r_1}u^{r_2}\ldots).
\end{align*}
Similarly,
$$
(u - i \tilde{u} )^{2k} = u^{2k} +(-i)^{2k}\tilde{u}^{2k} +
\sum_{(r_1, \ldots , r_{m})
\in \cal{S}}(-i)^{s(r)} (u^{r_1} \tilde{u}^{r_2}\ldots)+
 (-i)^{2k - s(r)} (\tilde{u}^{r_1}u^{r_2}\ldots).
$$
If $\cal{K} = \{r= (r_1,r_2,\ldots,r_{2m})\in
 \cal{S}; \ s(r) \in 2\N \}$, then
$$
(u + i \tilde{u})^{2k} + (u - i \tilde{u})^{2k} = 2u^{2k} +
 2(i)^{2k}\tilde{u}^{2k} + 2 \sum_{r\in \cal{K}}
 (i)^{s(r)} (u^{r_1}\tilde{u}^{r_2}\ldots) +
  (i)^{2k - s(r)} (\tilde{u}^{r_1}u^{r_2}\ldots),
$$
so from \eqref{3}, we get
$$
\Phi(d)^{2k}=
 \Phi(u^{2k}) + (i)^{2k}\Phi(\tilde{u}^{2k}) +
 \sum_{r \in \cal{K}}
(i)^{s(r)}  \Phi (u^{r_1} \tilde{u}^{r_2}\ldots) +
  (i)^{2k -s(r)} \Phi(\tilde{u}^{r_1}u^{r_2}\ldots).
$$

This implies
$$
(i)^{2k} \Phi(\tilde{u}^{2k}) = \Phi(d)^{2k} - \Phi(u^{2k}) -
 \sum_{r \in \cal{K}}(i)^{s(r)}
\Phi(u^{r_1}\tilde{u}^{r_2}\ldots) +
 i^{2k -s(r)} \Phi(\tilde{u}^{r_1}u^{r_2}\ldots).
$$
Taking the trace on both sides,
$$
|\T(\tilde{u}^{2k})| \leq |\T(d^{2k})| + |\T(u^{2k})| +
 \sum_{r  \in \cal{K}} | \T(u^{r_1} \tilde{u}^{r_2}\ldots)|+
 |\T(\tilde{u}^{r_1}u^{r_2}\ldots)|.
$$
Applying  Lemma~1,  with $\frac{1}{p_j} = \frac{r_j}{2k}$, for every
$r=(r_1, r_2, \ldots, r_{m}) \in \cal{K}$, we get
$$
|\T(\tilde{u}^{2k})| \leq |\T(d^{2k})| + \| u \|^{2k}_{2k} +
 \sum_{r  \in \cal{K}} (\| u \|^{r_1}_{2k} \| \tilde{u} \|^{r_2}_{2k}
\ldots) + (\| \tilde{u} \|^{r_1}_{2k} \| u \|^{r_2}_{2k}\ldots).
$$
We observe that by the definition of $\cal{K}$,
$$ \| u \|^{2k}_{2k} +
\sum_{r \in \cal{K}} (\| u \|^{r_1}_{2k} \|
\tilde{u} \|^{r_2}_{2k} \ldots) +
 ( \| \tilde{u}\|^{r_1}_{2k} \| u \|^{r_2}_{2k}\ldots)
$$
is  equal to  the sum of the terms of the expansion of $( \| u
\|_{2k} + \| \tilde{u} \|_{2k})^{2k}$ with $\| \tilde{u} \|_{2k}$ of
even exponents between $2$ and $2k-2$, i.e.,
\begin{multline*}
 \| u \|^{2k}_{2k} + \sum_{r
 \in \cal{K}} (\| u \|^{r_1}_{2k} \| \tilde{u} \|^{r_2}_{2k}
\ldots) + ( \| \tilde{u}\|^{r_1}_{2k}\| u \|^{r_2}_{2k} \ldots)\\
=\binom{2k}{0} \| u \|^{2k}_{2k} + \binom{2k}{2} \| u \|^{2k-2}_{2k}
\| \tilde{u} \|^{2}_{2k}
 + \cdots + \binom{2k}{2k-2} \| u \|^{2}_{2k} \| \tilde{u}
\|^{2k-2}_{2k}.
\end{multline*}
Since $\tilde{u}$ is self-adjoint, $\T(\tilde{u}^{2k})=
\| \tilde{u}\|_{2k}^{2k}$
and hence,
$$
\| \tilde{u} \|^{2k}_{2k} \leq \| d \|^{2k}_{2k} + \| u \|^{2k}_{2k} +
\binom{2k}{2} \| u \|^{2}_{2k} \| \tilde{u} \|^{2(k-1)}_{2k} + \cdots
+ \binom{2k}{2k-2} \| u \|^{2(k-1)}_{2k} \| \tilde{u} \|^{2}_{2k};
$$
and since $\| d \|^{2k}_{2k} \leq \| u \|^{2k}_{2k}$, we have
$$
\| \tilde{u} \|^{2k}_{2k} \leq 2 \| u \|^{2k}_{2k} + \binom{2k}{2} \|
u \|^2_{2k} \| \tilde{u} \|^{2(k-1)}_{2k} + \cdots + \binom{2k}{2k-2}
\| u \|^{2(k-1)}_{2k} \| u \|^2_{2k}.
$$

Divide both sides by $\|u \|^{2k}_{2k}$ and set $X_0 = \| \tilde{u} \|_{2k}/\|
u \|_{2k}$, we have
$$
X^{2k}_{0} - \binom{2k}{2} X^{2(k-1)}_{0} - \binom{2k}{4}
X^{2(k-2)}_{0} - \cdots - 2 \leq 0.
$$
Hence, $X_0$ is less than or equal to the largest real root of the
polynomial equation
$$
X^{2k} - \binom{2k}{2}X^{2(k-1)} - \binom{2k}{4}X^{2(k-2)} - \cdots -
2= 0.
$$
If the largest root is $K_{2k}$, we have
$$
\| \tilde{u} \|_{2k} \leq K_{2k} \| u \|_{2k}.
$$
 Using Minkowski's inequality, we conclude that for every $ f \in \cal{A}$
(not necessarily self adjoint), we have
$$
\| \tilde{f} \|_{2k} \leq 2K_{2k} \| f \|_{2k}.
$$
Since $\cal{A}$ is dense in $L^{2k}(\M,\T)$,
the inequality above shows that $\cal{H}$ can be
extended as a bounded linear operator from $L^{2k}(\M,\T)$ into
 $L^{2k}(\M,\T)$, so the theorem is proved for $p$ even.

For the general case, let
$ 2 \leq p < \infty$. Choose an integer $k$ such that $2k \leq p
\leq 2k + 2$. By \cite{DDP2} (Theorem~2.3), $L^p(\M,\T)$ can be
realized as a complex
interpolation of the pair of Banach spaces
 $(L^{2k}(\M,\T), L^{2k+2}(\M,\T))$, and we conclude
that $\cal{H}$ is also bounded from $L^p(\M,\T)$ into $L^p(\M,\T)$.

For $1 < p < 2$,  from the above case, $\cal{H}$ is bounded from $L^{q}(\M,\T)$
 into
$L^{q}(\M,\T)$, where $\frac{1}{p} + \frac{1}{q} = 1$, and we claim that as in
 the commutative case,  $(\cal{H})^* = - \cal{H}$.

To see this,
let  $u$ and  $v$ be self adjoint elements of   $ \cal{A}$; we have
$$
\Phi((u + i \tilde{u})(v + i \tilde{v})) = \Phi(u+i \tilde{u}) \cdot
\Phi(v + i \tilde{v}) = \Phi(u)\Phi(v),
$$
which implies that
$$
\Phi(uv + iu \tilde{v} + i \tilde{u} v - \tilde{u} \tilde{v}) =
\Phi(uv - \tilde{u} \tilde{v}) + i \Phi(u \tilde{v} + \tilde{u} v) =
\Phi(u)\Phi(v),
$$
so
$$
\T (uv - \tilde{u} \tilde{v}) + i \T(u \tilde{v} + \tilde{u} v) =
\T(\Phi(u)\Phi(v)).
$$
Since $\Phi(u)$ and $\Phi(v)$ are self-adjoint, $\T(\Phi(u) \Phi(v))
\in \R$, and also $\T (uv - \tilde{u} \tilde{v})$ and $\T(u \tilde{v}
+ \tilde{u} v) \in \R$. This implies $ \T(u \tilde{v} + \tilde{u}v)=0$
 and $ \T(u \tilde{v}) = - \T(\tilde{u} v)$. The proof is complete.
\qed

\begin{remarks}
 (1) One can deduce as in Corollary~2h of \cite{HR} that there exists a
 a constant $M$ such that $\|\tilde u\|_p \leq Mpq \|u\|_p$ for all
 $u \in L^p$, $1<p< \infty$, $1/p +1/q=1$.

 (2) We note that a result similar to that of Theorem~2 was also obtained by
Zsido (\cite{ZS}, Theorem~3.25) under different setting. The case of trace
class $c_p$ where $1< p <\infty$ was also considered by Asmar, Berkson and
Gillespie in \cite{ABG}. The proofs  given in \cite{ABG} and \cite{ZS} use
representations of locally compact groups. Our proof of Theorem~2 is more
direct.
 
 (3) After this paper was written, we learned that
Marsalli and West (Noncommutative $H^p$ spaces-Preprint)
also obtain Theorem~2 using similar approach.
\end{remarks}

\smallskip

Our next result can be viewed as a non-commutative variant of
Kolmogorov Theorem.
\begin{theorem}
Let  $u \in \M$ with $ u \geq 0$, and set $f = u + i \tilde{u}$.
Then for every
$s > 0$, $$ \T (\chi_{(s, \infty)} (|f|)) \leq 4 \frac{\| u
\|_1}{s}.$$
\end{theorem}

\smallskip
We will begin by collecting some lemmas necessary for the proof some of
which
are  probably known but some of the proofs will be included for the
convenience of the reader.

\begin{lemma}
For $u \in \M$, let $f = u +i\tilde{u}$ and $0<\varepsilon< 1$.
 The formal series
 $\displaystyle{\varphi(t)=\sum_{k=0}^\infty \frac{t^k}{k!}\varepsilon^k f^k} $
 is absolutely convergent in $H^1$ for every $t < 1/eM_1 \|u\|$ where
 $M_1=2M +2$ and $M$ is the constant from (1) of Remarks~1 above.
\end{lemma}
This lemma can be proved exactly as in Theorem~3a of \cite{HR}.
\begin{lemma}
For $u \in \M$, $u \geq 0$, let $f = u + i \tilde{u}$ and $ 0 < \varepsilon
< 1$.
\begin{itemize}
\item[(1)] $I + \varepsilon f$ has bounded inverse with $\|(I + \varepsilon
f)^{-1} \| \leq 1$.
\item[(2)] $f_{\varepsilon} = (\varepsilon I + f) (I + \varepsilon
f)^{-1} \in H^\infty$.
\item[(3)] $\Phi(f_\varepsilon) ={\Phi(u)}_\varepsilon $.
\item[(4)]$\Re\,( f_{\varepsilon}) \geq \varepsilon I$.
\item[(5)]$\Re\, (I +(f_\varepsilon -sI)(f_\varepsilon +sI)^{-1})\geq 0$
for every $s>0$.
\item[(6)] $\lim_{\varepsilon \to 0} \| f_{\varepsilon} -  f \|_p =0$
 \ ($ 1\leq p <\infty$).
\end{itemize}
\end{lemma}

\begin{proof}
(1) Note that $f$ is densely defined and that, for every $x \in
D(f)$,
$$
\langle(I + \varepsilon f) x, \, x \rangle = \langle(I+ \varepsilon
u)x, \, x \rangle + i \langle \tilde{u} x, \, x \rangle.
$$
Thus $|\langle(I+ \varepsilon f) x, \, x \rangle| \geq \|x\|^2$,
which implies
$
\|(I + \varepsilon f) x \| \geq \|x\| \quad \text{ for all } x \in
D(f)$.
So $I + \varepsilon f$ has bounded inverse with
 $\|(I + \varepsilon f)^{-1}\|\leq 1 $.
 \smallskip

(2) Note first that
$f_\varepsilon$ is bounded. In fact,
$f_{\varepsilon} = \varepsilon (I + \varepsilon f)^{-1} + f(I+
\varepsilon f)^{-1}$;
but $I = (I+\varepsilon f)(I+\varepsilon f)^{-1} = (I+ \varepsilon f)^{-1} +
\varepsilon f(I+\varepsilon f)^{-1}$ and $(I+\varepsilon f)^{-1}
\in \M$, so the operator $\varepsilon f(I+\varepsilon f)^{-1}$ is bounded,
implying $f(I+\varepsilon f)^{-1}$ is bounded.  To prove that
$f_{\varepsilon} \in H^\infty$, it suffices to show that
 $ (I +\varepsilon f)^{-1} \in H^\infty$.

Set $A=-\varepsilon f$. There exists a (unique) semi-group of contractions
$(T_t)_{t>0}$ such that $A$ is the infinitesimal generator of $(T_t)_{t>0}$
(see for instance \cite{Y} P.246-249). It is well known that
$$(I-A)^{-1}x= \int_{0}^\infty e^{-t} T_t x \ dt \quad \forall x \in H
\quad \text{and}$$
$$ T_t x= \lim_{n \to \infty}\exp\left( tA(I-n^{-1}A)^{-1}\right)x
\quad  \forall x \in H. $$

\noindent
Claim: $T_t \in H^\infty$ for every $t>0$.

Since $(T_t)_{t>0})$ is a semi-group, it is enough to verify this claim for
small $t$. Assume that $2t \leq 1/eM_1 \|u\|$. Let $\varphi(t)$ be the
operator in $H^1$ defined in Lemma~2. We will show that $T_t=\varphi(t)$.
Using the series expansion of the exponential and Lemma~2, we get
\begin{align*}
\|\exp\left(tA(I-n^{-1}A)^{-1}\right)-\phi(t)\|_1
&\leq \sum\limits_{k\geq 0} \frac{t^k}{k!}\|\left(A(I-n^{-1}A)^{-1}\right)^k
- A^k \|_1 \\
&\leq  \sum\limits_{k\geq 0}\| \frac{(2t)^k}{k!}\varepsilon^k
   f^k \|_1 <\infty.
\end{align*}
It is easy to verify that for each $k \geq 0$,
 $\lim\limits_{n \to \infty}
   \frac{t^k}{k!}\|\left(A(I-n^{-1}A)^{-1}\right)^k
- A^k \|_1  =0$ so
$$\lim_{n \to \infty}\|\exp\left(tA(I-n^{-1}A)^{-1}\right)-
\varphi(t)\|_1 =0 $$
which shows that $T_t =\varphi(t) \in H^1$ and the claim follows.

\noindent
We conclude the proof of (2) by noticing that $t \to T_t$  is
a continuous function in $H^1$ which shows that $(I-A)^{-1} \in H^1$.

\smallskip
(3) $f_\varepsilon  \in
H^\infty$ and $ f_\varepsilon \cdot (I + \varepsilon f) = \varepsilon I +
f \in H^2$,
so $ \Phi(f_{\varepsilon}) \Phi(I + \varepsilon f) = \Phi
(\varepsilon I + f)$. But $ \Phi(f) = \Phi(u)$, so we get $
\Phi(f_{\varepsilon}) (I + \varepsilon \Phi(u)) = \varepsilon I +
\Phi(u)$.

\smallskip
(4)\begin{align*}
 \Re\,(f_{\varepsilon}) &= \Re\,\left((\varepsilon I + f)(I +
\varepsilon f)^{-1} \right) \\
 &= \Re\, ([(\varepsilon I + \varepsilon^2 f) + (1
- \varepsilon^2) f)](I + \varepsilon f)^{-1}) \\
&= \varepsilon
I + (1 - \varepsilon^2) \Re\,(f(I + \varepsilon f)^{-1}).
\end{align*}
 Since we
assume that $\varepsilon < 1$, it is enough to show \noindent that $\Re\,(f(I+
\varepsilon f)^{-1}) \geq 0$. For this
\begin{align*}
\Re\, (f(I+\varepsilon f)^{-1}) &= \tfrac{1}{2}\left(f(I+\varepsilon
f)^{-1} + (I + \varepsilon f^*)^{-1}f^*\right) \\&= \tfrac{1}{2}(I +
\varepsilon f^*)^{-1} \left((I+\varepsilon f^*)f +
f^*(I+\varepsilon f)\right)(I+\varepsilon f)^{-1}  \\
&=\tfrac{1}{2}(I+\varepsilon f^*)^{-1}\left(2 \Re\, (f) + 2 \varepsilon
|f|^2\right)(I+\varepsilon f)^{-1} \geq 0.
\end{align*}

\smallskip
(5)
\begin{align*}
\Re\, (I &+ (f_{\varepsilon} - sI)(f_{\varepsilon} + sI)^{-1}) = I
+ \tfrac{1}{2} \left((f_{\varepsilon} - sI)(f_{\varepsilon} + sI)^{-1}
+ (f^*_\varepsilon + sI)^{-1}(f^*_\varepsilon - sI)\right) \\
&= I + \tfrac{1}{2}(f^*_\varepsilon + sI)^{-1} \left((f^*_\varepsilon
+ sI)(f_{\varepsilon} - sI) + (f^*_\varepsilon -
sI)(f_{\varepsilon} + sI)\right)(f_{\varepsilon} + sI)^{-1} \\
&= I + (f^*_\varepsilon + sI)^{-1} \left(|f_{\varepsilon}|^2 -
s^2I\right)(f_{\varepsilon} + sI)^{-1} \\
&= (f^*_\varepsilon + sI)^{-1} \left((f^*_\varepsilon +
sI)(f_{\varepsilon} + sI)+ |f_{\varepsilon}|^2 - s^2I\right)
(f_{\varepsilon} + sI)^{-1} \\
& = 2 (f^*_\varepsilon + sI)^{-1} \left(|f_{\varepsilon}|^2 + s
\Re\,f_{\varepsilon}\right)(f_{\varepsilon} + sI)^{-1},
\end{align*}
and the claim follows from the fact that $\Re\,(f_{\varepsilon}) \geq
\varepsilon I$.
\smallskip

(6)
We have for every $\varepsilon > 0$,
\begin{align*}
f_{\varepsilon} - f &= (\varepsilon I + f)(I +
\varepsilon f)^{-1} - f \\&= \left((\varepsilon I + f) - f(I+
\varepsilon f)\right)(I+\varepsilon f)^{-1} \\&= \varepsilon(I+
f^2)(I+\varepsilon f)^{-1},
\end{align*}
so
$$
\mu_t(f_\varepsilon - f) \pprec \varepsilon \mu_t(I+
f^2) \mu_t ((I+\varepsilon f)^{-1}).
$$
Since $\| (I + \varepsilon f)^{-1}\| \leq 1$, we get
$\mu_t((I+ \varepsilon f)^{-1}) \leq 1$ for every $t>0$.
 Also $I + f^2 \in
L^p(\M,\T)$ for every $p > 1$, so $\| f_{\varepsilon} - f \|_p
\leq \varepsilon \| I + f^2 \|_p \to 0$ (as $ \varepsilon \to 0$).
The proof is complete.
\end{proof}

\begin{lemma}
Let  $a$ and $b$ be operators in
 $\overline{ \M}$ with $ a \geq 0$, $b \geq 0$, and let $P$ be a
projection that commutes with $a$. Then $\T(ab) \geq \T(P(ab)P)$.
\end{lemma}
\begin{proof}
To see this, notice that, since $P$ commutes with $a$, $PaP \leq a$,
so $b^{1/2}PaPb^{1/2} \leq b^{1/2}ab^{1/2}$, implying that
$\T(b^{1/2}PaPb^{1/2} ) \leq \T(b^{1/2}ab^{1/2})$ and
\begin{align*}
\T(P(ab)P) &= \T(P(ab)) \\
&= \T(PaPb) = \T(b^{1/2}PaPb^{1/2}) \\
&\leq \T(b^{1/2}a b^{1/2}) = \T(ab).
\end{align*}
The lemma is proved.
\end{proof}

\begin{lemma}
Let $S$ be a positive operator that commutes with $|f_{\varepsilon}|$ then
$$\T(S \Re\, (f_\varepsilon))\leq \T(S|f_\varepsilon|).$$
\end{lemma}

\begin{proof} Let $S\geq 0$ and assume that
$S|f_\varepsilon|=|f_\varepsilon|S$. we have
\begin{align*} \T(S\Re\,(f_\varepsilon))
 &=\T(S(f_\varepsilon +f_{\varepsilon}^*)/2) \\
 &=1/2(\T(Sf_\varepsilon) +\overline{\T(Sf_\varepsilon)}) \\
 &\leq |\T(Sf_\varepsilon)|
\end{align*}
Let $f_\varepsilon =u|f_\varepsilon|$  be the polar decomposition of
$f_\varepsilon$. Since $S$ commutes with $|f_\varepsilon|$, we get
$|\T(Sf_\varepsilon)|=|\T(u|f_\varepsilon|S)|
=|\T(uS^{1/2}|f_\varepsilon|S^{1/2})|\leq \T(S|f_\varepsilon|)$. Thus the
proof is complete.
\end{proof}

\begin{lemma}
Let $A$ and $B$ be positive operators such that:
\begin{itemize}
\item[(i)] $A^{-1}$ and $B^{-1}$ exists;
\item[(ii)] $\T(CA) \leq \T(CB)$ for every $C$ that commutes with $B$;
\end{itemize}
Then for every $C$ that commutes with $B$, $\T(C B^{-1}) \leq \T(C A^{-1})$
\end{lemma}

\begin{proof}
Observe that $\T (C B^{-1}) = \T(\alpha \beta )$,
where $\alpha = C^{1/2} B^{-1} A^{1/2}$ and $ \beta = A^{-1/2}
C^{1/2}$. By H\"older's inequality, %

\begin{align*}
\T (C B^{-1}) &\leq \T( |\alpha|^2)^{1/2} \, \T(|\beta|^2)^{1/2} =
\T(\alpha^* \alpha)^{1/2}\, \T(\beta^* \beta)^{1/2} = \T(\alpha
\alpha^*)^{1/2}  \T(\beta^* \beta)^{1/2} \\
&= \T(C^{1/2} B^{-1} A^{1/2} A^{1/2} B^{-1} C^{1/2} )^{1/2}
\, \T(C^{1/2} A^{-1/2} A^{-1/2} C^{1/2} )^{1/2} \\
&= \T(C^{1/2} B^{-1} A
B^{-1} C^{1/2})^{1/2}\,  \T(C A^{-1} )^{1/2}\\
&=\T(B^{-1}CB^{-1}A)^{1/2}\, \T(CA^{-1})^{1/2}.
\end{align*}
Since $C$ commutes with $B$, the operator $B^{-1}CB^{-1}$ commutes with $B$
so
 we get by assumption that $\T(B^{-1}C B^{-1} (A)) \leq
 \T(B^{-1}C B^{-1} (B))$,
and therefore
 $$\T(C B^{-1}) \leq \T(C B^{-1})^{1/2}\, \T(C
A^{-1})^{1/2}$$
 which shows that  $\T(C B^{-1}) \leq \T (C A^{-1}).$
The proof of the lemma is complete.
\end{proof}

\smallskip
\noindent
{\bf Proof of Theorem~3.}

 Our proof is inspired by  the argument of Helson in \cite{HEL}
for the commutative case.

Let $u$ and $f$ be as in the statement of the Theorem~3, and
fix  $0 < \varepsilon < 1$. Set $f_{\varepsilon} $ as in  Lemma~3.
For $ s \in (0, \infty)$ fixed, consider the following transformation
on $\{z; \Re\,(z) \geq 0 \}$:
$$
A_s(z) = 1 + \frac{z-s}{z+s} \quad \text{ for all } z \in \{ w, \Re\,(w)
\geq 0 \}.
$$
It can be checked that the part of the plane $\{ z; |z| \geq s
\}$ is mapped to the half disk $\{ w;\,\Re\,(w) \geq 1 \}$;
this fact is very crusial in the argument of \cite{HEL} for the
 commutative case.

Note that $\sigma (f_{\varepsilon} )$ is a compact subset of $\{ z;
\Re\,(z) \geq \varepsilon \}$. By the analytic
 functional calculus for Banach algebras, %
$$
A_s (f_{\varepsilon}) = I + (f_{\varepsilon} - s
I)(f_{\varepsilon} +  sI)^{-1} \in H^\infty
$$
and therefore (since $A_s$ is analytic)
\begin{equation}
\eqlabel{*}
\Phi(A_s(f_{\varepsilon})) = A_s
(\Phi(f_{\varepsilon}))=A_s(\Phi(u)_\varepsilon).
\end{equation}

Note that since $\Phi(u)$ is self-adjoint, so are
$\Phi(u)_\varepsilon$ and $A_s(\Phi(u)_\varepsilon)$.
 We conclude  from \eqref{*}
that $\T\left(I+ (\Phi(u)_\varepsilon - sI)
(\Phi(u)_\varepsilon + sI)^{-1}\right) \in \R$, and therefore
\begin{equation}\eqlabel{1}
\T\left(
\Re\,(I +(f_{\varepsilon} - sI)(f_{\varepsilon} + sI)^{-1})\right) =
\T(A_s(\Phi(u)_\varepsilon)).
\end{equation}

Let $P=\chi_{(s,\infty)}(|f_{\varepsilon}|)$.
 The projection $P$ commutes with $|f_{\varepsilon}|$ and we  have
$$
\Re \left[I + (f_{\varepsilon} - sI)(f_{\varepsilon} + sI)^{-1}\right] =
(f^*_\varepsilon + sI)^{-1} \left[2|f_{\varepsilon}|^2 + 2 s
\Re\, (f_{\varepsilon})\right] (f_{\varepsilon} + sI)^{-1};
$$
but since $\Re\,(f_{\varepsilon}) \geq \varepsilon I \geq 0$, we get
$$
2|f_{\varepsilon}|^2 + 2s \Re\,(f_{\varepsilon}) \geq
2|f_{\varepsilon}|^2 ,
$$
and hence %
\begin{equation}
\eqlabel{**}
\T\left[\Re\,(I + (f_{\varepsilon} - sI)(f_{\varepsilon}+sI)^{-1})\right] \geq
\T \left[2|f_{\varepsilon}|^2(f_{\varepsilon} + sI)^{-1}
(f_{\varepsilon}^* + sI)^{-1}\right].
\end{equation}

Applying Lemma~4 for $a = 2|f_{\varepsilon}| $ and
$b = (f_{\varepsilon} + sI)^{-1} (f^*_\varepsilon + s
I)^{-1}$, we obtain %
\begin{equation}
\eqlabel{***}
\T\left[\Re\,(I + (f_{\varepsilon} -sI) (f_{\varepsilon} + sI)^{-1})\right]
\geq \T [ 2P |f_{\varepsilon}|^2
(f_{\varepsilon} + sI)^{-1} (f^*_\varepsilon + sI)^{-1}].
\end{equation}

Note that $(f_{\varepsilon} + sI)^{-1}(f^*_\varepsilon +
sI)^{-1} = \left(|f_{\varepsilon}|^2 + 2s \Re\,(f_{\varepsilon}) +
s^2I \right)^{-1}$.

Set
\begin{align*}
A &= (f^*_\varepsilon + sI) (f_{\varepsilon} + sI)
= |f_{\varepsilon}|^2 + 2s \Re\,(f_{\varepsilon}) +s^2I
, \\
B &= |f_{\varepsilon}|^2 + 2s |f_{\varepsilon}|+ s^2I.
\end{align*}
It is easy to see from Lemma~5 that if $C$ is a
positive  operator that commutes
with $B$ then $\T(CA)\leq \T(CB)$.

 Applying Lemma~6 to
$A$, $B$ and $C = 2P|f_{\varepsilon}|^2 $,
we obtain ( from \eqref{***})  that
$$
\T \left[\Re\, (I + (f_{\varepsilon} - sI) (f_{\varepsilon} + sI)^{-1})\right]
\geq \T (C A^{-1}) \geq \T (C B^{-1})
$$
and hence
\begin{equation}
\eqlabel{****}
\T \left[ \Re\, (I + (f_{\varepsilon} - sI)(f_{\varepsilon} + sI)^{-1})\right]
\geq \T \left[2P|f_{\varepsilon}|^2
(|f_{\varepsilon}|^2 + s|f_{\varepsilon}| + s^2I)^{-1}\right].
\end{equation}

If we denote by $E^{|f_{\varepsilon}|}$ the spectral decomposition
of $|f_{\varepsilon}|$, then
$$
2P|f_{\varepsilon}|^2
(|f_{\varepsilon}|^2 + 2 s |f_{\varepsilon}| + s^2 I )^{-1}
= \int^{\infty}_s \frac{2t^2 }{t^2 + 2st + s^2}\, d
E^{|f_{\varepsilon}|}_t.
$$
Let
$$\psi_{s}(t) = \frac{2t^2}{t^2 + 2 st
+ s^2} \quad \text{ for } t \in [s, \infty).
$$
One can show that $\psi_{s} $ is increasing on  $[s, \infty)$ so
$\psi_{s}(t)\geq \psi_{s}(s)=1/2$  for $t\geq s$,
and therefore
$$
\int^\infty_s \frac{2t^2 }{t^2 + 2 s t + s^2}\, d
E^{|f_{\varepsilon}|}_t \geq \frac{1}{2}P,
$$
so we deduce from \eqref{****} that
$$
\T\left[\Re\,(I + (f_{\varepsilon} - sI)(f_{\varepsilon} +
sI)^{-1})\right]\geq \frac{1}{2}\T(P).
$$

To finish the proof, recall from \eqref{1} that
$$
\T\left[\Re\,(I+(f_{\varepsilon}-sI)(f_{\varepsilon}+sI)^{-1})\right]=
\T\left[I+(\Phi(u)_\varepsilon - sI)(\Phi(u)_\varepsilon + sI)^{-1}\right],
$$
so
\begin{align*}
\T(P) &\leq 2
 \T\left[I + (\Phi(u)_\varepsilon
- sI) (\Phi(u)_\varepsilon + sI)^{-1}\right] \\
&= 2 \T\left[2 \Phi(u)_\varepsilon
(\Phi(u)_\varepsilon + sI)^{-1}\right].
\end{align*}
But $(\Phi(u)_\varepsilon + sI)^{-1} =
 \frac{1}{s} (\frac{\Phi(u)_\varepsilon}{s} +
I)^{-1} $ has norm $\leq 1/s$, hence%
$$
\T(P) \leq 4 \frac{\| u_\varepsilon \|_1
}{s}.
$$
Now taking  $\varepsilon \to 0$,
 we get from Lemma~3~(6) that
  $\| u_\varepsilon \|_1 \to \| u\|_1$ and
  $\| f_\varepsilon -f \|_1 $ converges to zero. In particular,
   $(f_\varepsilon)$ converges to $f$ in measure. We obtain from
\cite{FK}(Lemma~3.4) that
 $\mu_t(f)\leq \liminf\limits_{n \to \infty}\mu_t(f_{\varepsilon_n})$ for
 each $t>0$ and $\varepsilon_n \to 0$. This implies that for every
 $s>0$ and every $t>0$,
  $\chi_{(s, \infty)}\left(\mu_t(f)\right) \leq
  \liminf\limits_{n \to \infty}\chi_{(s, \infty)}
  \left(\mu_t(f_{\varepsilon_n})\right)$. Hence by Fatou's lemma,
  \begin{align*}
\T(\chi_{(s, \infty)}(|f|))&=
\int_{0}^1 \chi_{(s, \infty)}\left(\mu_t(|f|)\right)\ dt \\
&\leq \liminf\limits_{n \to \infty}
\int_{0}^1 \chi_{(s, \infty)}\left(\mu_t(|f_{\varepsilon_n}|)\right)\ dt \\
&=\liminf\limits_{n \to \infty}\T(\chi_{(s, \infty)}(|f_{\varepsilon_n}|))\\
&\leq \limsup\limits_{n \to \infty}
(4||u_{\varepsilon_n} ||_1/s).
\end{align*}

Hence $\T(\chi_{(s, \infty)} (|f|)) \leq 4\, {\| u \|_1 }/s$.
The proof is complete. \qed

\bigskip
We are now ready to extend the Hilbert transform to $L^1(\M,\T)$.
Recall that $L^{1, \infty}(\M,\T) = \{ a \in \overline{\M};\  \sup_{t > 0} t
\mu_t (a) < \infty \}$.

 Set $\|a\|_{1, \infty } = \sup_{t>0}t
\mu_t(a)$ for $a \in L^{1, \infty}(\M, \T)$. As in the commutative case,
$\| . \|_{1, \infty}$ is equivalent to a quasinorm in $L^{1, \infty}
(\M, \T)$, so there is a fixed constant $C$ such that, for every $a, b
\in L^{1, \infty} (\M, \T)$, we have $\|a + b \|_{1, \infty} \leq C(\| a
\|_{1, \infty} + \| b \|_{1, \infty})$.

For $u \in \M$, let $Tu = u + i \tilde{u}$. From Theorem 1, $T$ is
linear and Theorem~3 can be restated as follows:
$$
\text{ For any } u \in \M \text{ with } u \geq 0, \text{ we have } \|
Tu \|_{1, \infty} \leq 4 \| u \|_1;
$$
this implies that for $u \geq 0$,
$$
\| \tilde{u} \|_{1, \infty} \leq C ( 4+1) \| u \|_1 = 5C \| u \|_1.
$$
Now suppose  that $u \in \M$, $u = u^*$, $u = u_+ - u_-$ and $\tilde{u}
= \tilde{u}_+ - \tilde{u}_-$. Then
$$
\| \tilde{u} \|_{1, \infty} \leq C( \| \tilde{u}_+ \|_{1,\infty} +
\| \tilde{u}_- \|_{1,\infty} ) \leq 5C^2 \| u \|_1.
$$
Similarly, if we require only $u \in \M$, %
we have $u = \Re\,(u)+ i \Im\, (u)$ and by linearity,
 $\tilde{u} = \widetilde{\Re(u) } + i\,
\widetilde{\Im(u)}$, and as above,
$$
\| \tilde{u} \|_{1, \infty} \leq 10C^3 \| u \|_1.
$$
We are now ready to define the extension $\cal{H}$ in $L^1(\M,\T)$: If $u \in
L^1(\M,\T)$, let $(u_n)_{n \in  \N}$ be a sequence in $\M$ such that $\| u -
u_n \|_1 \to 0$ as $n \to \infty$. Then
$$
\| \tilde{u}_n - \tilde{u}_m \|_{1, \infty} \leq 10 C^3 \| u_n - u_m \|_1,
$$
and since $\| u_n - u_m \|_{1} \to 0$ as $n,m \to \infty$, the
sequence $(\tilde{u}_n)_n$ converges in $L^{1, \infty}(\M, \T)$ to an
operator $\tilde{u}$. This defines $\tilde{u}$ for $u \in L^1(\M,\T) $.
 This
definition can be easily checked to be independent of the sequence
$(\tilde{u}_n)_n$ and agree with the conjugation operator  defined for
$p >1 $.

 Letting $n \to \infty$ in the inequality $\| \tilde{u}_n
\|_{1, \infty} \leq 10 C^3 \| u_n \|_1$, we obtain the following
theorem ($H^{1,\infty}$ denotes the closure of $H^\infty$ in
 $L^{1,\infty}(\M,\T)$):
\begin{theorem}
There is a unique extension of $\cal{H}$ from $L^1(\M,\T) $
into $L^{1, \infty}(\M,\T)$ with the following property: $u + i \tilde{u}
\in H^{1, \infty}$ for all $ u \in L^1(\M,\T)$, and there is a constant
$K$ such that $\| \tilde{u} \|_{1, \infty } \leq K \| u \|_1$ for all
$ u \in L^1(\M,\T)$.
\end{theorem}

\smallskip
\begin{corollary}
For any $p$ with $0 < p < 1$ there exists a constant $K_p$ such that
$$\|\tilde{u}\|_p \leq K_p \|u\|_1 \ \ \text{for all}\  u \in L^1(\M,\T).$$
\end{corollary}

\begin{proof}
It is enough to show that such a constant exists for $u \in \M$, $u \geq 0$.
Recall that for $u \in \M$, the distribution $\lambda_s(u)$ equals
$\T(\chi_{(s, \infty)} (u))$.

 Let $F(s) = 1 - \lambda_s(u) = \T
(\chi_{(0,s)}(u))$. Assume that $\|u\|_1 \leq 1$. From Theorem~3,
$$
1- F(s) \leq \frac{4}{s}\|u\|_1 \leq \frac{4}{s}.
$$
Note that $F$ is a non-increasing right continuous function and for $p
> 0$, %
$$
\T (|\tilde{u}|^p) = \int^1_0 \mu_t(|\tilde{u}|)^p d t =
\int^\infty_0 s^p d F(s) \leq 1 + \int^\infty_1 s^p d F(s).
$$
If $A$ is
a point of continuity for $F ( A > 1)$, then
$$
\int^A_1 s^p d F(s) = [s^p(F(s) - 1)]^A_1 + p\int^A_1 (1- F(s))s^{p-1}ds.
$$
Since $1 - F(s) \leq \frac{4}{s}$, we get that both $[s^p(F(s)
-1)]^A_1$ and $\int^A_1 (1 - F(s))s^{p-1} d s$ are bounded for $0 < p
< 1$, that is, $\int^1_0 \mu_t(|\tilde{u}|)^p dt$ has bound
independent of $u$.
\end{proof}

\smallskip
The {\bf Riesz projection} $\cal{R}$ can now be defined as
 in the commutative case:
for every $ a \in L^p(\M,\T)$, \ $(1\leq p \leq \infty)$,
   $$\cal{R}(a) = \frac{1}{2}(a + i\tilde{a} + \Phi(a)).$$

From Theorem~2, one can easily verify that $\cal{R}$ is
 a bounded projection from $L^p(\M,\T)$ onto $H^p$ for $1<p<\infty$.
 In particular $H^p$ is a complemented subspace of $L^p(\M,\T)$.

  For
$p=1$, Theorem~4 shows that $\cal{R}$ is bounded from $L^1(\M,\T)$  into
$H^{1,\infty}$.

Our next result gives a sufficient condition on an operator $a\in L^1(\M,\T)$
so that its conjugate $\tilde{a}$ belongs to $L^1(\M,\T)$.
\begin{theorem}
 There exists a constant $K$ such that  for every 
positive $a \in \M$,
 $$ \| \tilde{a}\|_1 \leq K\T( a\log^+a) + K.$$
\end{theorem}

\begin{proof}
Let  $C$ be the absolute consrant such that
 $\| \tilde{a}\|_p \leq Cpq \|a\|_p$ for all
$a \in L^p(\M,\T)$, $1<p<\infty$ and $1/p +1/q=1$.
 The conclusion of the theorem
can be deduced as a straightforward adjustment of the commutative case
in \cite{ZYG}(vol~{.II}, p.~{119}); we will present it here for completeness.

  Let $a \in \M$; we will assume first
 that $a\geq 0$. Let $(e_t)_t $ be the spectral
decomposition of $a$. For each $k \in \N$, let $P_k =\chi_{[2^{k-1},2^k)}(a)$
be the  spectral projection relative to $[2^{k-1}, 2^k)$. Define $a_k =aP_k$
 for
$k \geq 1$ and $a_0=a\chi_{[0, 1)}(a)$.
 Clearly
$a= \sum_{k=0}^{\infty} a_k $ in $L^p(\M,\T)$ for every $1<p<\infty$.
 By linearity, $\tilde{a}=\sum_{k=0}^\infty \tilde{a}_k$. For every $k\in\N$,
 $\|\tilde{a}_k\|_1 \leq \|\tilde{a}_k\|_p \leq Cp^2{(p-1)}^{-1} \|a_k\|_p$.
 Since $a_k \leq 2^k P_k$, we get for $1<p<2$,
 $$\|\tilde{a}_k\|_1 \leq 4C \frac{1}{p-1}2^k \T(P_k)^{\frac{1}{p}}.$$
If we set $p =1 +\frac{1}{k+1}$ and $\epsilon_k =\T(P_k)$, we have
 $$ \|\tilde{a}_k\|_1 \leq 4C (k+1)2^k \epsilon_{k}^{\frac{k+1}{k+2}}.$$
Taking the summation over $k$,
 $$\|\tilde{a}\|_1 \leq \sum_{k=0}^\infty 4C(k+1)2^k \epsilon_{k}^{\frac{k+1}{k+2}}.$$
We note as in \cite{ZYG} that if $J=\{ k \in \N;\ \epsilon_k \leq 3^{-k}\}$
then
 $$\sum_{k \in J} 4C(k+1)2^k \epsilon_{k}^{\frac{k+1}{k+2}} \leq
 \sum_{k=0}^\infty 4C(k+1)2^k (3^{-k})^{\frac{k+1}{k+2}} = \alpha <\infty.$$
On  the other hand, for $k \in \N \setminus J$,
$\epsilon_{k}^{\frac{k+1}{k+2}} \leq \epsilon_k. 3^{\frac{k}{k+2}} \leq
\beta \epsilon_k$ where $\beta=\sup_k 3^{\frac{k}{k+2}}$. So we get
\begin{align*}
\|\tilde{a}\|_1 &\leq
 \alpha + 4C\beta \sum_{k=0}^\infty (k+1)2^k \epsilon_k \\
 &\leq \alpha +4C\beta(\epsilon_0 + 4\epsilon_1) +4C\beta \sum_{k\geq 2}(k+1)
2^k \epsilon_k.
\end{align*}
Since for $k\geq 2$,\   $k+1 \leq 3(k+1)$, we get
$$\|\tilde{a}\|_1 \leq
\alpha +16C\beta +24C\beta\sum_{k\geq2}(k-1)2^{k-1}\epsilon_k.$$
To complete the proof, notice that  for $k\geq 2$,
\begin{align*}
 (k-1)2^{k-1} \epsilon_k &= \int_{2^{k-1}}^{2^k}  (k-1)2^{k-1}\ d\T(e_t)\\
  &\leq \int_{2^{k-1}}^{2^k} \frac{t\log t}{\log 2} \ d\T(e_t).
\end{align*}
Hence by setting $K=\max\{\alpha +16C\beta, 24C\beta/\log 2\}$, we get:
 $$\|\tilde{a}\|_1 \leq K + K \tau\left(a \log^+(a)\right).$$

 The proof is complete.
\end{proof}
\begin{remark}
From Theorem~4, one can deduce that if 
$a \in L^1(\M,\T)$ is such that $a\geq 0$ and
$a\log^+(a)$ belongs to $L^1(\M,\T)$ then $\tilde{a} \in L^1(\M,\T)$.
\end{remark}

\smallskip
Let us finish with the following open question that arises naturally from
the commutative case and the topic of this paper:
it is a well known result of Bourgain (\cite{BO6}) that
 $L^1(\mathbb{T})/H^1(\mathbb{T})$ is of cotype 2 and later
 Lancien (\cite{LAN}) proved that a similar result holds for $L^1/H^1$
 associated with weak*-Dirichlet algebras.

\smallskip
\noindent
{\bf Problem:}\  Is $L^1(\M,\T)/H^1(\M,\T)$ of cotype 2? (or merely of
finite cotype?)

It should be noted that the theory of conjugate functions and the
boundedness of the Riesz projection were very crusial in the proof given by
Bourgain (\cite{BO6}) for the classical case
 and Lancien (\cite{LAN}) for the
setting of weak*-Dirichlet algebras.
\smallskip

\noindent
{\bf Acknowledgements.}
Parts of the work reported here were done during the author's visit
at the Mathematical Sciences Research Institute (MSRI) in Berkeley, CA.
The author is grateful to the organizers of the Special Semester on
Convex Geometry for financial support.


\end{document}